\documentclass{article}
\usepackage{graphicx} 

\usepackage{amsmath}
\usepackage{amssymb}
\usepackage{amsthm}
\usepackage[colorlinks=true]{hyperref}

\title{An Explicit Uniform Mordell Conjeture over Function Fields of Characteristic Zero}
\author{Jiawei Yu}
\date{\today}

\newtheorem{thm}{Theorem}[section]
\newtheorem{cor}[thm]{Corollary}
\newtheorem{lem}[thm]{Lemma}
\newtheorem{prop}[thm]{Proposition}

\begin{document}

\maketitle

\tableofcontents

\section{Introduction}

Let $C$ be a geometrically connected smooth projective curve of genus $g>1$ over $\mathbb{Q}$. Mordell \cite{Mordell} conjectured that $C(\mathbb{Q})$ is finite. This conjecture is proved by Faltings \cite{Fal} with $\mathbb{Q}$ replaced by any number field. Before Faltings' proof, the analog of this conjecture for curves over complex function fields was proved independently by Manin \cite{Manin} and Grauert \cite{Grauert}. Vojta \cite{Vojta_Geometric,Vojta_Ar} gave an alternative proof with Diophantine approximation in both cases. Based on Vojta's proof, Dimitrov-Gao-Habegger \cite{DGH} and Kühne \cite{Kuhne} proved the following uniform theorem (cf. \cite[Theorem 1.1]{Gao}).

\begin{thm}[Dimitrov-Gao-Habegger+Kühne]
For any integer $g>1$, there is a constant $c(g)$ with the following property. Let $K$ be a field of characteristic $0$, $C/K$ a geometrically connected smooth projective curve of genus $g$, $J$ the Jacobian variety of $C/K$, and $P_0\in C(K)$ a rational point. Then for any subgroup $\Gamma\subseteq J(K)$ of finite rank $\rho$,
$$\sharp(i_{P_0}(C(K))\cap\Gamma)\le c(g)^{\rho+1}.$$
\end{thm}

Here $i_{P_0}$ is the Abel-Jacobi map
$$i_{P_0}:C\longrightarrow J,\quad P\longmapsto P-P_0.$$
The rank of $\Gamma$ is the dimension of the $\mathbb{Q}$-vector space $\Gamma\otimes\mathbb{Q}$. Note that $\Gamma$ is not required to be finitely generated. In this article, we determine the constant $c(g)$ in the non-isotrivial case explicitly. Our main theorem is as follows.

\begin{thm}
\label{Main}
Let $K$ be a field of characteristic $0$, $C/K$ a geometrically connected smooth projective curve of genus $g>1$, $J$ the Jacobian variety of $C/K$, and $\alpha$ a line bundle on $C$ of degree $1$. If $C$ is non-isotrivial over $\mathbb Q$, then for any subgroup $\Gamma\subseteq J(K)$ of finite rank $\rho$,
$$\sharp(i_\alpha(C(K))\cap\Gamma)\le(16g^2+32g+184)(20g)^{\rho}.$$
In particular, if $J(K)$ is of finite rank $\rho$, then
$$\sharp C(K)\le(16g^2+32g+184)(20g)^{\rho}.$$
\end{thm}

A curve $C$ is \emph{non-isotrivial} over $\mathbb{Q}$ if it is not isomorphic to the base change of a curve from $\bar{\mathbb{Q}}$ to $\bar{K}$. We denote
$$i_\alpha:C\longrightarrow J,\quad P\longmapsto P-\alpha.$$
Then $i_{P_0}=i_{\mathcal{O}(P_0)}$ is a special case. In our proof, there is a natural $\alpha$ which need not to be a point in $C$.

We will reduce the main theorem to the following theorem.

\begin{thm}
\label{main}
Let $k$ be an algebraically closed field of characteristic $0$, $K$ the function field of a smooth projective connected curve $B/k$, $C/K$ a smooth projective geometrically connected curve non-isotrivial over $k$ and of genus $g>1$, $J$ the Jacobian variety of $C/K$, and $\alpha$ a line bundle on $C$ of degree $1$. Then for any subgroup $\Gamma\subseteq J(\bar K)$ of rank $\rho$,
$$\sharp(i_\alpha(C(\bar K))\cap\Gamma)\le(16g^2+32g+184)(20g)^{\rho}.$$
In particular, if the $K/k$-trace of $J$ is trivial, then
$$\sharp C(K)\le(16g^2+32g+184)(20g)^{\rho_{LN}}.$$
Here the Lang-Néron rank $\rho_{LN}$ is the rank of $J(K)$.
\end{thm}

The \emph{$K/k$-trace} $(\mathrm{tr}_{K/k}(J),\mathrm{tr})$ of an abelian variety $A/K$ is the final object in the category of pairs $(B,f)$ consisting of an abelian variety $B/k$ and a morphism $f:B_K\to A$. If $k$ is algebraically closed in $K$, the $K/k$-trace exists and $J(K)/\mathrm{tr}_{K/k}(J)(k)$ is finitely genereated  (cf. \cite[Theorem 2.1 and 6.2]{Conrad}). In particular, our assumption on trivial trace implies that $\rho_{NL}$ is finite.

With the curve $B$, we have a Weil height on $C(\bar K)$. The rational points are divided into two parts based on height and counted separately. Note that we divide points according to the admissible pairing of $C$ introduced in \cite{Zhang_Admissible} rather than Faltings height, which was commonly used (cf. \cite{Gao}). Yuan \cite{Yuan_Admissible} showed that these two quantities on curves bound each other uniformly in the function field case.

We modify Vojta's proof of the Mordell conjecture to count points with large height. Vojta constructed an effective divisor with the Grothendieck-Riemann-Roch theorem. Then he gave a Diophantine approximation inequality on its index to derive an ineffective upper bound of height. An obstruction to giving explicit inequality is that the Weil height can be determined only up to a bounded function. Zhang \cite{Zhang_adelic} introduced the adelic line bundle, the height associated to which is the canonical height. We use it to refine Vojta's estimates and bound uniformly the number of points with large height. We also use Siu's theorem \cite{Siu} to avoid dealing with higher cohomologies and simplify the construction of the divisor. 

On the other hand, a conjecture of Bogomolov states that there are only finitely many points in $C(\bar{K})$ with small height. Ullmo \cite{Ullmo} proved it with the equidistribution theorem of Szpiro-Ullmo-Zhang \cite{SUZ}. Zhang \cite{Zhang_Admissible} showed that the strict positivity of the admissible pairing implies the Bogomolov conjecture. Based on it, Looper-Silverman-Wilms \cite{LSW} gave a quantitative uniform result on the Bogomolov conjecture over function fields, and Yuan \cite{Yuan_Admissible} proved uniform Bogomolov conjecture over global fields independently. We combine the theorem of Looper-Silverman-Wilms and the modified Vojta's inequality to deduce our main theorem.

Throughout this article, let $\mathcal{C}$ be the minimal regular model of $C$ over $B$. Since replacing $K$ by a finite extension does not change the first conclusion of Theorem \ref{main}, we may assume that $C$ has split semi-stable reduction. The dualizing sheaf $\bar{\omega}=\omega_{\mathcal{C}/B}$ is a line bundle on $\mathcal{C}$. Let $X=C\times_KC$ be the product, $p_1,p_2:X\to C$ the projections, and $\Delta\subseteq X$ the diagonal.

We are grateful to Xinyi Yuan for his vital suggestions and patient help. We thank Zheng Xiao and Chunhui Liu for teaching the author Diophantine approximation. We thank Yinchong Song for reviewing the draft of this article. We thank Joseph Silverman and Robert Wilms for helpful comments. We thank Ziyang Gao for helpful conversation. We thank the anonymous referees for helpful suggestions.

\section{Difference between heights}

In this section we recall the adelic line bundles introduced in \cite{Zhang_adelic} and admissible metrics in \cite{Zhang_Admissible}. We refer to \cite[Appendix]{Yuan_Admissible} for a detailed treatment on admissible metrics. Then we compare the height associated to $\bar{\omega}$ with the canonical height.

For $v\in B$, let $K_v$ be the local field and $\mathcal{O}_{K_v}$ its valuation ring. We take $\vert\varpi_v\vert=e^{-1}$ for simplicity. Here $\varpi_v$ is a uniformizer of $K_v$. For a projective variety $Z/K$ and a line bundle $L$ on $Z$, a model $(\mathcal{Z},\mathcal{L})$ of $(Z,L^{\otimes n})$ over $\mathcal{O}_{K_v}$ induces a metric $\Vert\cdot\Vert$ on $L_{K_v}$ as follows: For $z\in Z(\bar{K}_v)$, it extends to
$$\bar{z}:\mathrm{Spec}(\mathcal{O}_{\bar{K}_v})\longrightarrow\mathcal{Z}.$$
For $\ell\in z^*L$, define
$$\Vert\ell\Vert=\inf_{a\in\bar{K}_v}\{\vert a\vert:\ell^n\in a^n\bar{z}^*\mathcal{L}\}.$$
A metric $\Vert\cdot\Vert'$ on $L_{K_v}$ is continuous and bounded if $\Vert\cdot\Vert'/\Vert\cdot\Vert$ is continuous and bounded for some metric $\Vert\cdot\Vert$ induced by a model. An \emph{adelic metric} on $L$ is a collection $\{\Vert\cdot\Vert_v\}$ of continuous and bounded metrics $\Vert\cdot\Vert_v$ on $L_{K_v}$ for all $v\in B$, such that $\Vert\cdot\Vert_v$ is induced by a model of $(Z,L)$ over an open subvariety $U\subseteq B$ for all $v\in U$. An \emph{adelic line bundle} is a pair $\bar{L}=(L,\{\Vert\cdot\Vert_v\})$ consisting of a line bundle $L$ and an adelic metric $\{\Vert\cdot\Vert_v\}$ on $L$.

We say that an adelic metric  $\{\Vert\cdot\Vert_{v}\}$ on $L$ is the limit of a sequence of adelic metrics $\{\Vert\cdot\Vert_{n,v}\}(n=1,2,\dots)$, if $\Vert\cdot\Vert_{n,v}$ is independent of $n$ for $v$ in some open subvariety of $B$, and $\Vert\cdot\Vert_{n,v}/\Vert\cdot\Vert_{v}$ convenges uniformly on $X(\bar{K_v})$ for each $v\in B$. A model $(\mathcal{Z},\mathcal{L})$ of $(Z,L^{\otimes{n}})$ on $B$ is relatively nef if $\mathcal{L}$ is nef on special fibers of $\mathcal{Z}$. An adelic metric $\{\Vert\cdot\Vert_v\}$ is relatively nef if it is the limit of a sequence of adelic metrics induced by relatively nef models over $B$. An adelic line bundle is \emph{integrable} if it is the tensor quotient of two relatively nef adelic line bundles. Denote by $\widehat{\mathrm{Pic}}(Z)_{\mathrm{int}}$ the group of isometry classes of integrable adelic line bundles on $Z$.

If $Z=\mathrm{Spec}(K)$, the \emph{degree} of $\bar{L}\in\widehat{\mathrm{Pic}}(Z)_{\mathrm{int}}$ is defined as
$$\widehat{\deg}(\bar{L})=\sum_{v\in B}-\log\Vert s\Vert_v,$$
where $s$ is any non-zero section of $L$. If $d=\dim{Z}$, there is a \emph{Deligne pairing} 
$$\widehat{\mathrm{Pic}}(Z)_{\mathrm{int}}^{d+1}\longrightarrow\widehat{\mathrm{Pic}}(K)_{\mathrm{int}},\quad(\bar{L}_1,\dots,\bar{L}_{d+1})\longmapsto\pi_*\langle\bar{L}_1,\dots,\bar{L}_{d+1}\rangle$$
with respect to the structure morphism $\pi:Z\to\mathrm{Spec}(K)$ (cf. \cite{YZ_Adelic}). The intersection number of $d+1$ line bundles is the degree of their Deligne pairing. The \emph{height} associated to $\bar{L}$ is 
$$h_{\bar{L}}:Z(\bar K)\longrightarrow\mathbb{R},\quad z\longmapsto\frac{\widehat{\deg}(z^*\bar{L})}{\deg(z)}.$$
Here $\deg(z)$ is the degree of the residue field of $z$ over $K$.

The space $\hat{H}^0(\bar{L})$ of \emph{effective sections} of an adelic line bundle $\bar{L}=(L,\{\Vert\cdot\Vert_v\})$ on $Z$ is consists of section $s$ of $L$ such that
$$\sup\Vert s\Vert_v\le1$$
for all $v$. When $\bar{L}$ is induced from $\mathcal{L}$ on some projective model, $\hat{H}^0(\bar{L})=H^0(\mathcal{L})$. The \emph{volume} of $\bar{L}$ is
$$\mathrm{vol}(\bar{L})=\lim_{n\to\infty}\frac{(d+1)!}{n^{d+1}}\dim_k\hat{H}^0(\bar{L})$$
for $d=\dim(Z)$. The limit always exists. If $\mathrm{vol}(\bar{L})>0$, then we say $\bar{L}$ is \emph{big}. For nef adelic line bundles $\bar{L}_1$ and $\bar{L}_2$, there is Siu's inequality
$$\mathrm{vol}(\bar{L}_1-\bar{L}_2)\ge\bar{L}_1^{d+1}-(d+1)\bar{L}_1^d\bar{L}_2.$$

Choose a line bundle $\alpha_0$ of degree $1$ on $C$ satisfying $(2g-2)\alpha_0=\omega_{C/K}$. Let $\theta$ be the image of
$$C^{g-1}\longrightarrow J,\quad (x_1,\dots,x_{g-1})\longmapsto i_{\alpha_0}(x_1)+\dots+i_{\alpha_0}(x_{g-1}).$$
It is a divisor on $J$. The line bundle $\Theta=\mathcal{O}(\theta)+[-1]^*\mathcal{O}(\theta)$ is symmetric. Zhang \cite{Zhang_adelic} applied Tate's limiting argument to construct an integrable adelic metric on $\Theta$. Denote the adelic line bundle by $\bar{\Theta}$. The \emph{canonical height} is defined as the associated height
$$\hat{h}=h_{\bar{\Theta}}:J(\bar K)\longrightarrow \mathbb{R}.$$
It is positive and $\vert\cdot\vert=\hat{h}(\cdot)^{1/2}$ extends to a norm on
$$J(\bar K)_{\mathbb{R}}=J(\bar K)\otimes\mathbb{R}$$
satisfying the parallelogram law. Denote the corresponding inner product by $\langle\cdot,\cdot\rangle_{\Theta}$. By abuse of language, we write $\hat{h}(i_{\alpha_0}(x))$ as $\hat{h}(x)$.

For $v\in B$, let $\Gamma_v$ be the reduction graph of the special fiber $\mathcal{C}_v$, i.e. the vertexes and edges of $\Gamma_v$ represent the components and nodes of $\mathcal{C}_v$ respectively, and each edge is of length $1$. Denote by $F(\Gamma_v)$ the space of continuous and piecewise smooth function on $\Gamma_v$. For $f\in F(\Gamma_v)$, the Laplacian operator gives a measure
$$\Delta f=-f''(x)dx-\sum d_{\overrightarrow v}f(P)\delta_P.$$
Here $x$ represents a canonical coordinate on each edge. The summation is over $P\in\Gamma_v$ and tangent directions $\overrightarrow v$ at $P$, and $\delta_P$ is the Dirac measure supported at $P$.

For each probability measure $\mu$, there exists a unique symmetric function $g_\mu:\Gamma_v^2\to\mathbb{R}$, called the \emph{Green function} associated to $\mu$, satisfying
$$g_\mu(x,\cdot)\in F(\Gamma),$$
$$\Delta g_\mu(x,\cdot)=\delta_x-\mu,$$
$$\int_{\Gamma_v}g_\mu(x,\cdot)\mu=0.$$
It can be seen directly from definition that
$$g_\mu(x,x)=\sup_{y\in\Gamma_v}g_\mu(x,y)\ge 0.$$

The \emph{canonical divisor} $K_{\Gamma_v}$ of $\Gamma_v$ is a formal linear combination
$$K_{\Gamma_v}=\sum \deg(\bar\omega\vert_{F_\xi})\xi.$$
The summation is over all vertexes $\xi$, and $F_\xi$ is the component of $\mathcal{C}_v$ represented by $\xi$. There is a unique probability measure $\mu$ satisfying $g_\mu(K_{\Gamma_v},x)+g_\mu(x,x)$ is a constant independent of $x$. Denote the measure by $\mu_v$ and the Green function by $g_v$. With the retraction map $C(\bar{K_v})\to\Gamma_v$ (cf. \cite[Chapter 4]{Ber}), we can view $g_v$ as a function on $C(\bar{K_v})^2$.

The model $(\mathcal{C},\bar{\omega})$ induces an adelic metric $\{\Vert\cdot\Vert_{\mathrm{Ar},v}\}$ on $\omega_{C/K}$. The \emph{canonical admissible metric} $\{\Vert\cdot\Vert_{a,v}\}$ of $\omega_{C/K}$ defined by
$$\Vert\cdot\Vert_{a,v}(x)=\Vert\cdot\Vert_{\mathrm{Ar},v}(x)\cdot e^{g_v(x,x)},\quad x\in C_{K_v}^{\mathrm{an}}.$$
is an integrable adelic metric. Denote by $\omega_a$ the canonical admissible adelic line bundle. The self-intersection number $\omega_a^2$ is non-negative.

Similarily, there is an integrable canonical admissible adelic line bundle
$$\mathcal{O}(\Delta)_a=(\mathcal{O}(\Delta),\{\Vert\cdot\Vert_{\Delta,v}\})\in\widehat{\mathrm{Pic}}(X)_{\mathrm{int}}$$
determined by
$$\Vert\cdot\Vert_{\Delta,a}(x,y)=e^{-i_v(x,y)-g_v(x,y)},\quad x,y\in C(\bar{K_v}),\,x\ne y.$$
Here, $i_v(x,y)$ is the stable intersection number, i.e. if $x,y\in C(K_w')$ for some finite extension $K_w'/K_v$, and $\bar{x},\bar{y}$ are their closures in the minimal regular model $\mathcal{C'}$ over $\mathcal{O}_{K_w'}$, then
$$i_v(x,y)=\frac{(\bar{x}\cdot\bar{y})}{[K_w':K_v]}.$$

The following proposition gives the difference we need. Note that the global stable intersection number
$$i(P_1,P_2)=\sum_{v\in B}i_v(P_1,P_2)\log(e)=\frac{(\bar{P_1}\cdot\bar{P_2})}{[K':K]}$$
can be considered as a Weil height function associated to the line bunlde $\mathcal{O}(\Delta)$ at the point $P=(P_1,P_2)\in (X\backslash\Delta)(\bar{K})$. So we can view the second equality as the Mumford's equality in Vojta's proof.

\begin{prop}
\begin{itemize}
    \item [(1)] {For $P_0\in C(\bar{K})$,
$$h_{\bar\omega}(P_0)=\frac{g-1}{g}\hat{h}(P_0)+\sum_{v\in B}g_v(P_0,P_0)+\frac{\omega_a^2}{4g(g-1)}.$$}
    \item [(2)] {For $P_1,P_2\in C(\bar{K})$, if $P_1\ne P_2$, then
$$i(P_1,P_2)=\frac{\hat{h}(P_1)}{2g}+\frac{\hat{h}(P_2)}{2g}-\langle P_1,P_2\rangle_{\Theta}-\sum_{v\in B}g_v(P_1,P_2)-\frac{\omega_a^2}{4g(g-1)}.$$}
\end{itemize}
\end{prop}

\begin{proof}
(1) Consider the morphism
$$i_\omega:C\longrightarrow J,\quad x\longmapsto(2g-2)x-\omega_{C/K}.$$
By \cite[Theorem 2.10(1)]{Yuan_Admissible},
$$i_\omega^*\bar{\Theta}=4g(g-1)\omega_a-\pi^*\pi_*\langle\omega_a,\omega_a\rangle$$
in $\widehat{\mathrm{Pic}}(C)\otimes\mathbb Q$. Here $\pi_*\langle\cdot,\cdot\rangle$ is the Deligne pairing for the structure morphism $\pi:C\to\mathrm{Spec}(K).$
Since $i_\omega=[2g-2]\circ i_{\alpha_0}$, we have
$$(2g-2)^2\hat{h}(P_0)=4g(g-1)h_{\omega_a}(P_0)-\omega_a^2.$$
Hence,
\begin{align*}
h_{\bar\omega}(P_0)=&\sum_{v\in B}g_v(P_0,P_0)+h_{\omega_a}(P_0)\\
=&\sum_{v\in B}g_v(P_0,P_0)+\frac{g-1}{g}\hat{h}(P_0)+\frac{\omega_a^2}{4g(g-1)}.
\end{align*}

(2) Consider the morphism
$$j:X\longrightarrow J,\quad(x,y)\longmapsto y-x.$$
By \cite[Theorem 2.10(2)]{Yuan_Admissible},
$$j^*\bar{\Theta}=2\mathcal{O}(\Delta)_a+p_1^*\omega_a+p_2^*\omega_a$$
in $\widehat{\mathrm{Pic}}(X)\otimes\mathbb Q$. Hence
$$\hat{h}(P_2-P_1)=2h_{\mathcal{O}(\Delta)_a}(P_1,P_2)+h_{\omega_a}(P_1)+h_{\omega_a}(P_2).$$
We get
\begin{align*}
i(P_1,P_2)=&h_{\mathcal{O}(\Delta)_a}(P_1,P_2)-\sum_{v\in B}g_v(P_1,P_2)\\
=&\sum_{i=1,2}\frac{\hat{h}(P_i)-h_{\omega_a}(P_i)}{2}-\langle P_1,P_2\rangle_{\Theta}-\sum_{v\in B}g_v(P_1,P_2)\\
=&\frac{\hat{h}(P_1)}{2g}+\frac{\hat{h}(P_2)}{2g}-\frac{\omega_a^2}{4g(g-1)}-\langle P_1,P_2\rangle_{\Theta}-\sum_{v\in B}g_v(P_1,P_2).
\end{align*}
\end{proof}

Denote by $\delta(\Gamma_v)$ the total length of $\Gamma_v$. Zhang \cite{Zhang_phi} introduced the $\varphi$-invariant as
$$\varphi(\Gamma_v)=-\frac{1}{4}\delta(\Gamma_v)+\frac{1}{4}\int_{\Gamma_v}g_v(x,x)((10g+2)\mu_v-\delta_{K_{\Gamma_v}}).$$
He also showed in the loc. cit. that
$$\omega_a^2\ge\frac{2g-2}{2g+1}\sum_{v\in B}\varphi(\Gamma_v).$$

By \cite[Lemma 2.2]{LSW} and \cite[Proposition 13.7]{BR}, we have
$$\vert g_v(x,y)\vert\le\frac{15}{4}\varphi(\Gamma_v).$$
Therefore,
$$\vert\sum_{v\in B}g_v(x,y)\vert\le\frac{30g+15}{8g-8}\omega_a^2$$
As a consequence, we have the following corollary.

\begin{cor}
\label{difference}
\begin{itemize}
    \item [(1)]{For $P_0\in C(\bar{K})$,
$$0\le h_{\bar\omega}(P_0)-\frac{g-1}{g}\hat{h}(P_0)\le\frac{19}{2}\omega_a^2.$$}
    \item[(2)]{For $P_1,P_2\in C(\bar{K})$, if $P_1\ne P_2$, then
$$-\frac{19}{2}\omega_a^2\le i(P_1,P_2)-\frac{\hat{h}(P_1)}{2g}-\frac{\hat{h}(P_2)}{2g}+\langle P_1,P_2\rangle_{\Theta}\le\frac{37}{4}\omega_a^2.$$}
\end{itemize}
\end{cor}

\section{Vojta's inequality}
In this section we modify Vojta's proof \cite{Vojta_Geometric} to give a uniform inequality. The following is the main theorem of this section.

\begin{thm}
\label{ineq}
For $P_1,P_2\in C(\bar K)$, if 
$$\vert P_2\vert\ge20g\vert P_1\vert$$
and
$$\vert P_1\vert\ge200\sqrt{g\omega_a^2},$$
then
$$\frac{\langle P_1,P_2\rangle_\Theta}{\vert P_1\vert\vert P_2\vert}\le\frac{4}{5}.$$
\end{thm}

Let $\bar{M}$ be the pull-back to $X$ of an adelic line bundle on $B$ of degree $1$. Consider the adelic line bundle
$$\bar{L}=d_1p_1^*{\omega_a}+d_2p_2^*{\omega_a}+d((2g-2)\mathcal{O}({\Delta})_a-p_1^*{\omega_a}-p_2^*{\omega_a})+c\bar{M},$$
where $d_1,d_2,d,c$ are positive integers to be decided.

\begin{lem}
If $d_1\ge gd$, $d\ge d_2$, and $2(d_1d_2-gd^2)c>(g+1)d_1d^2\omega_a^2$, then $L$ is big.
\end{lem}

\begin{proof}
Take
$$\bar{L}_1=(d_1-d)p_1^*{\omega_a}+d(g-1)p_2^*{\omega_a}+d(2g-2)\mathcal{O}({\Delta})_a+c\bar{M},$$
$$\bar{L}_2=(gd-d_2)p_2^*{\omega_a}.$$
By \cite[Theorem 2.10(2)]{Yuan_Admissible},
$$p_1^*\omega_a+p_2^*\omega_a+2\mathcal{O}(\Delta)_a$$
is nef. So both $\bar{L}_1$ and $\bar{L}_2$ are nef. By Siu's inequality, we have
\begin{align*}
&\mathrm{vol}(\bar{L})\\
\ge&\bar{L}_1^3-3\bar{L}_1^2\bar{L}_2\\
=&6(2g-2)^2(d_1d_2-gd^2)c+d^3(2g-2)^2\left((2g+1)\omega_a^2-(2g-2)\sum_{v\in B}\phi(\Gamma_v)\right)+\\
&3(2g-2)(d_1^2d_2-(g^2+2g-1)d_1d^2+2gd^3+(2g-2)d_1dd_2-(2g-1)d^2d_2){\omega_a^2}\\
\ge&6(2g-2)^2(d_1d_2-gd^2)c-3(2g-2)(g^2+g)d_1d^2{\omega_a^2}\\
>&0.
\end{align*}
Here we use the intersection numbers in the proof of \cite[Theorem 3.6]{Yuan_Admissible}
\end{proof}


For $P_i\in C(K)(i=1,2)$, we have $P=(P_1,P_2)\in X(K)$ and sections $\bar{P_i}\subset\mathcal{C}$, $\bar{P}\subset\mathcal{X}$ of $B$. Choose a local coordinate $x_i$ on $C$ at $P_i$. Then for any effective divisor $D\in\mathrm{Div}(X)$, $D$ is defined near $P$ by a formal power series
$$\sum_{i_1,i_2\ge0}a_{i_1,i_2}x_1^{i_1}x_2^{i_2}.$$
Recall that the \emph{index} of $D$ at $P$ with respect to a pair of positive numbers $(e_1,e_2)$ is
$$\mathrm{ind}(D,P,e_1,e_2)=\mathrm{min}\{\frac{i_1}{e_1}+\frac{i_2}{e_2}:a_{i_1,i_2}\ne 0\}.$$
It is independent of the choice of $x_1,x_2$.

\begin{proof}[Proof of Theorem \ref{ineq}.]
The theorem is invariant after replacing $K$ by a finite extension. We may assume $P_1,P_2\in C(K)$. Take
$$d_1=\lceil\sqrt{g+\frac{1}{400}}\frac{\vert P_2\vert}{\vert P_1\vert}d\rceil,$$
$$d_2=\lceil\sqrt{g+\frac{1}{400}}\frac{\vert P_1\vert}{\vert P_2\vert}d\rceil,$$
$$c=\lceil200(g+1)d_1\omega_a^2\rceil,$$
where $\lceil\cdot\rceil$ is the ceiling function, i.e. $\lceil x\rceil$ is the least integer not less than $x$. For $d$ large enough, $\bar{L}$ is big. There is a positive integer $n$ such that $n\bar{L}$ admits an effective section $s$. 

Let $\mathcal{X}$ be the regular locus of $\mathcal{C}\times_B\mathcal{C}$. For $v\in B$, there is a natural bijection $(\mathcal{C}\times_B\mathcal{C})(\mathcal{O}_{\bar K_v})=X(\bar K_v)$. Denote by $X(\bar K_v)^\circ$ the subset of $X(\bar{K}_v)$ corresponding to $\mathcal{X}(\mathcal{O}_{\bar K_v})$. It is an analytic domain of $X(\bar{K}_v)$. Let $\bar{\Delta}$ be the closure of $\Delta$ in $\mathcal{X}$ and
$$\mathcal{L}=d_1p_1^*\bar{\omega}+d_2p_2^*\bar{\omega}-d((2g-2)\mathcal{O}(\bar{\Delta})-p_1^*\bar{\omega}-p_2^*\bar{\omega})+cM\in\mathrm{Pic}(\mathcal{X}).$$
Then $\mathcal{L}_K=L$. For $v\in B$, $\mathcal{L}$ induces a metric $\Vert\cdot\Vert_{\mathcal{L},v}$ of $L$ on $X(\bar K_v)^\circ$ as in the projective case. Let
$$m_v=-\log(\sup_{x\in X(\bar K_v)^\circ}\Vert s(x)\Vert_{n\mathcal{L},v}).$$
View $s$ as a rational section of $n\mathcal{L}$. Then $m_v$ is the minimal multiplicity of the vertical part of $\mathrm{div}(s)$ over $v$. So
$$D=\mathrm{div}(s)-\sum_{v\in B}m_v[\mathcal{X}_v]$$
is an effective divisor on $\mathcal{X}$. Here $\mathcal{X}_v$ is the fiber of $\mathcal{X}\to B$. Denote
$$\mathcal{L}'=n\mathcal{L}-\sum_{v\in B}m_v\mathcal{O}([\mathcal{X}_v]).$$
Since $s$ is an effective section of $\bar{L}$,
$$m_v\ge-\log(\sup\Vert s\Vert_{n\bar{L},v})-n(d_1+d_2+2gd)\sup|g_v|\ge-2nd_1\sup |g_v|.$$
So
$$\sum_{v\in B}m_v\ge-20nd_1\omega_a^2\ge-nc.$$



By \cite[Lemma 4.1]{Vojta_Geometric}, there is a finite extension $B'$ of $B$ with function field $K'$ and a regular surface $\mathcal{C}_i'$ over $B'$ for $i=1,2$ satisfying

\begin{itemize}
    \item[(1)]{there is a morphism $\mathcal{C}_i'\to\mathcal{C}\times_BB'$;}
    \item[(2)]{its restriction to the generic fiber $C_i'\to C\times_KK'$ is of degree $2d_{3-i}$ and unramified outside $P_i$;}
    \item[(3)]{there are two points of $C_i'$ lying over $P_i$, both defined over $K'$ and of ramification index $d_{3-i}$.}
\end{itemize}
Let $\mathcal{X}'$ be the regular locus of $\mathcal{C}_1'\times_{B'}\mathcal{C}_2'$. Then there is a morphism
$$f:\mathcal{X}'\longrightarrow\mathcal{X}\times_BB'$$
extending $C_1\times_{K'}C_2\to X\times_KK'$. Choose arbitary $P'=(P_1',P_2')\in C_1'\times_{K'}C_2'$ lying over $P$. The closure $\bar{P}'$ of $P'$ in $\mathcal{X}'$ is proper. The conormal sheaf ${\mathcal{N}}_{\bar{P}'/\mathcal{X}'}^\vee$ is a direct sum
$${\mathcal{N}}_{\bar{P}'/\mathcal{X}'}^\vee={\mathcal{N}}_{\bar{P}_1'/\mathcal{C}_1'}^\vee\oplus{\mathcal{N}}_{\bar{P}_2'/\mathcal{C}_2'}^\vee.$$
By \cite[Lemma 3.2]{Vojta_Dyson},
\begin{align*}
&\mathrm{ind}(D,P,(2g-2)nd_1,(2g-2)nd_2)\\
=&\mathrm{ind}(f^*D,P',(2g-2)nd_1d_2,(2g-2)nd_1d_2)\\
=&\frac{\mathrm{ind}(f^*D,P',1,1)}{(2g-2)nd_1d_2}.
\end{align*}
Together with \cite[Lemma 4.2.2]{Vojta_Geometric}, we have
\begin{align*}
&\mathrm{ind}(D,P,(2g-2)nd_1,(2g-2)nd_2)\\
\ge&\frac{1}{(2g-2)nd_1d_2}\frac{-2\deg(f^*\mathcal{L}'\vert_{\bar P'})}{\deg(\omega_{\mathcal{C}_1'/B'}\vert_{\bar P_1'})+\deg(\omega_{\mathcal{C}_2'/B'}\vert_{\bar P_2'})}\\
=&\frac{1}{(g-1)n}\frac{-\deg(\mathcal{L}'\vert_{\bar P})}{d_1\deg(\bar{\omega}\vert_{\bar P_1})+d_2\deg(\bar{\omega}\vert_{\bar P_2})}\\
\ge&\frac{1}{g-1}\frac{-\deg(\mathcal{L}\vert_{\bar P})-c}{d_1\deg(\bar{\omega}\vert_{\bar P_1})+d_2\deg(\bar{\omega}\vert_{\bar P_2})}\\
=&-\frac{1}{g-1}+\frac{1}{g-1}\frac{d(\deg(\bar{\omega}\vert_{\bar P_1})+\deg(\bar{\omega}\vert_{\bar P_2})-(2g-2)\deg(\mathcal{O}(\tilde{\Delta})\vert_{\bar P}))-2c}{d_1\deg(\bar{\omega}\vert_{\bar P_1})+d_2\deg(\bar{\omega}\vert_{\bar P_2})}.
\end{align*}
By Corollary \ref{difference},
$$\frac{g-1}{g}\vert P_i\vert^2\le\deg(\bar{\omega}\vert_{\bar P_i})\le\frac{g-1}{g}\vert P_i\vert^2+\frac{19}{2}\omega_a^2\le\frac{100}{99}\frac{g-1}{g}\vert P_i\vert^2.$$
Note that $\deg(\mathcal{O}(\tilde{\Delta})\vert_{\bar P})=i(P_1,P_2)$. Again by Corollary \ref{difference},
$$\deg(\bar{\omega}\vert_{\bar P_1})+\deg(\bar{\omega}\vert_{\bar P_2})-(2g-2)\deg(\mathcal{O}(\tilde{\Delta})\vert_{\bar P})\ge(2g-2)(\langle P_1,P_2\rangle_{\Theta}-\frac{19}{2}\omega_a^2).$$
We may assume $\langle P_1,P_2\rangle_\Theta\ge0$. Then
\begin{align*}
&\mathrm{ind}(D,P,(2g-2)nd_1,(2g-2)nd_2)\\
\ge&-\frac{1}{g-1}+\frac{1}{g-1}\frac{(2g-2)d\langle P_1,P_2\rangle_{\Theta}-(19g-19)d\omega_a^2-2c}{d_1\deg(\bar{\omega}\vert_{\bar P_1})+d_2\deg(\bar{\omega}\vert_{\bar P_2})}\\
\ge&-\frac{1}{g-1}+\frac{0.99g}{g-1}\frac{2d\langle P_1,P_2\rangle_\Theta}{d_1\vert P_1\vert^2+d_2\vert P_2\vert^2}-\frac{g}{(g-1)^2}\frac{(19g-19)d\omega_a^2+2c}{d_1\vert P_1\vert^2+d_2\vert P_2\vert^2}\\
\ge&-\frac{1.01}{g-1}+\frac{0.99g}{g-1}\frac{2d\langle P_1,P_2\rangle_\Theta}{d_1\vert P_1\vert^2+d_2\vert P_2\vert^2}-\frac{g}{(g-1)^2}\frac{2c}{d_1\vert P_1\vert^2+d_2\vert P_2\vert^2}.
\end{align*}
As a consequence,
\begin{align*}
&\liminf_{d\to\infty}\,\mathrm{ind}(D,P,(2g-2)nd_1,(2g-2)nd_2)\\
\ge&-\frac{1.01}{g-1}+\frac{0.99g}{(g-1)\sqrt{g+\frac{1}{400}}}\frac{\langle P_1,P_2\rangle_{\Theta}}{\vert P_1\vert\vert P_2\vert}-\frac{g}{(g-1)^2}\frac{200(g+1)\omega_a^2}{\vert P_1\vert^2}\\
\ge&-\frac{1.03}{g-1}+\frac{0.98g}{(g-1)\sqrt{g}}\frac{\langle P_1,P_2\rangle_{\Theta}}{\vert P_1\vert\vert P_2\vert}.
\end{align*}

By Dyson's lemma (cf. \cite{Vojta_Dyson}),
$$V(\mathrm{ind}(D,P,(2g-2)nd_1,(2g-2)nd_2))\le\frac{d_1d_2-gd^2}{d_1d_2}+(2g-1)\frac{d_2}{2d_1}.$$
Here
$$V(t)=\int_{x,y\in[0,1],x+y\le t}dxdy.$$
Taking limit we have
\begin{align*}
&\liminf_{d\to\infty}\,V(\mathrm{ind}(D,P,(2g-2)nd_1,(2g-2)nd_2))\\
\le&\frac{1}{400g+1}+(2g-1)\frac{\vert P_1\vert^2}{2\vert P_2\vert^2}\\
\le&\frac{1}{200g}
\end{align*}
By the monotonicity of $V$, we have
$$\liminf_{d\to\infty}\,\mathrm{ind}(D,P,(2g-2)nd_1,(2g-2)nd_2)\le\frac{1}{10\sqrt{g}}.$$

Combining two inequalities on index we have
$$-\frac{1.03}{g-1}+\frac{0.98g}{(g-1)\sqrt{g}}\frac{\langle P_1,P_2\rangle_{\Theta}}{\vert P_1\vert\vert P_2\vert}\le\frac{1}{10\sqrt{g}}.$$
Therefore,
$$\frac{\langle P_1,P_2\rangle_{\Theta}}{\vert P_1\vert\vert P_2\vert}\le\frac{100}{98}(\frac{g-1}{10g}+\frac{1.03}{\sqrt{g}})\le\frac{4}{5}.$$
\end{proof}

\section{Proof of the main theorem}

In this section we finish the proof of Theorem \ref{Main} and Theorem \ref{main}. 

\begin{proof}[Proof of Theorem \ref{Main}]
Since $\Gamma$ is of finite rank, there is a finitely generated subgroup $\Gamma_0\subseteq\Gamma$ satisfying
$$\Gamma_0\otimes\mathbb{Q}=\Gamma\otimes\mathbb{Q}.$$
Then we can find a finitely generated extension $K/\mathbb{Q}$ such that $C$ and $\Gamma_0$ are defined over $K$. For any $P\in\Gamma$, there is $Q\in\Gamma_0$ and an integer $n$ such that $nP=Q$. The morphism $[n]:J\mapsto J$ is finite. Therefore $P\in J(\bar K)$.

Denote by $M_g$ the coarse moduli scheme of curves of genus $g$ over $\bar{\mathbb{Q}}$. Let
$$\iota:\mathrm{Spec}(\bar K)\longrightarrow M_g$$
be the $\bar K$-point corresponding to $C_{\bar{K}}$. Since $C$ is non-isotrivial over $\mathbb Q$, $x$ does not factor throught $\mathrm{Spec}{(\bar{\mathbb{Q}})}$. Note that
$$\cap k=\bar{\mathbb{Q}},$$
where the intersection is over all algebraically closed subfield $k\subseteq\bar K$ satisfying $\bar K/k$ is of transcendental degree $1$. There is a $k$ such that $\iota$ does not factor through $\mathrm{Spec}(k)$. Replace $K$ by $Kk$. Then $K/k$ is finitely generated of transcendental degree $1$, and hence the function field of a smooth projective connected curve $B/k$. Theorem \ref{Main} follows from Theorem \ref{main}.
\end{proof}

We need the following proposition to count points with large height in the proof of Theorem \ref{main}

\begin{prop}
\label{gap}
Let $P_1,P_2\in C(\bar K)$ be two distinct points. If
$$\vert P_2\vert\ge\vert P_1\vert\ge\left(5\sqrt{2}-\frac{1}{4}\right)\sqrt{\omega_a^2}$$
and
$$\frac{\langle P_1,P_2\rangle_{\Theta}}{\vert P_1\vert\vert P_2\vert}\ge\frac{4}{5},$$
then
$$\frac{\vert P_2\vert}{\vert P_1\vert}\ge\frac{9g}{10}.$$
\end{prop}

\begin{proof}
Since $P_1\ne P_2$, by Corollary \ref{difference} we have
\begin{align*}
0&\le i(P_1,P_2)\\
&\le\frac{\vert P_1\vert^2}{2g}+\frac{\vert P_2\vert^2}{2g}-\langle P_1,P_2\rangle_{\Theta}+\frac{37}{4}\omega_a^2\\
&\le\frac{\vert P_1\vert^2}{2g}+\frac{\vert P_2\vert^2}{2g}-\frac{4}{5}\vert P_1\vert\vert P_2\vert+\frac{37}{4}\omega_a^2\\
&\le\frac{\vert P_1\vert^2}{2g}+\frac{\vert P_2\vert^2}{2g}-\frac{3}{5}\vert P_1\vert\vert P_2\vert.
\end{align*}
Solving it we have
$$\frac{\vert P_2\vert}{\vert P_1\vert}\ge\frac{3g}{5}+\sqrt{\left(\frac{3g}{5}\right)^2-1}\ge\frac{9g}{10}.$$
\end{proof}

\begin{proof}[Proof of Theorem \ref{main}]
Denote by $\Gamma_1$ the subgroup of $J(\bar{K})$ generated by $\Gamma$ and $\alpha-\alpha_0$. We have the vector spaces $V=\Gamma\otimes\mathbb{R}$ and $V_1=\Gamma_1\otimes\mathbb{R}$. Let 
$$W=V+\alpha-\alpha_0$$
be the coset of $V$ in $V_1$. Then
\begin{align*}
\sharp(i_{\alpha}(C(\bar K))\cap\Gamma)=&\sharp\{P\in C(\bar K):P-\alpha\in V\}\\
=&\sharp\{P\in C(\bar K):P-\alpha_0\in W\}.
\end{align*}

For any point $x\in V_1$ and positive number $r$, denote by $B(x,r)$ the closed ball with center $x$ and radius $r$ in $V_1$. By \cite{LSW}, a ball of radius
$$R_0=\sqrt{\frac{\omega_a^2}{8(g^2-1)}}$$
contains at most there are at most $16g^2+32g+124$ rational points in $C(\bar{K})$. Let
$$R_1=\left(5\sqrt{2}-\frac{1}{4}\right)\sqrt{\omega_a^2}.$$
We cover $B(O,R_1)\cap W$ with balls of radius $R_0$ in an inductive manner. Choose arbitary $x_1\in B(O,R_1)\cap W$. After selecting $x_1,\dots,x_t$, if $B(O,R_1)\cap W$ is not covered by $B(x_1,R_1),\dots,B(x_t,R_0)$, choose arbitary
$$x_{t+1}\in (B(O,R_1)\cap W)-(\bigcup_{i=1}^tB(x_i,R_0)).$$
Note that $B(x_i,R_0/2)\cap W$ are disjoint and all contained in $B(O,R_1+R_0/2)\cap W$. Since $W$ is the translation of a $\rho$-dimensional vector space. Calculating volume in $W$, we find that $B(O,R_1)\cap W$ can be covered by at most
$$\frac{(R+R_0/2)^\rho}{(R_0/2)^\rho}\le(20g)^\rho$$
balls of radius $R_0$. Therefore,
$$\sharp\{P\in C(\bar K):P-\alpha_0\in V, \vert P\vert\le R_1\}\le(16g^2+32g+124)(20g)^\rho.$$

The above inequality implies the first assertion in Theorem \ref{main} in the case $\rho=0$. Now suppose $\rho>0$. For any $0\le\phi\le\pi/2$, by \cite[Lemma 6.3]{Vojta_Geometric}, $V_1$ can be covered by
$$\frac{2\rho}{\sin(\phi/4)^{\dim(V_1)}\cos(\phi/4)}$$
sectors such that any two points in the same sector form an angle of at most $\phi$. Here $\dim(V_1)\le\rho+1$. Let $\phi=\arccos(4/5)$. Then
$$\frac{2\rho}{\sin(\phi/4)^{\dim(V_1)}\cos(\phi/4)}\le13^{\rho+1}.$$
Since
$$\frac{200\sqrt{g\omega_a^2}}{R_1}<\left(\frac{9g}{10}\right)^7,$$
by Proposition \ref{gap}, for any sector $S$,
$$\sharp\{P\in C(\bar K):P-\alpha\in S,R_1\le\vert P\vert\le200\sqrt{g\omega_a^2}\}\le7.$$
By Theorem \ref{ineq}, similarly we have
$$\sharp\{P\in C(\bar K):P-\alpha\in S,\vert P\vert\ge200\sqrt{g\omega_a^2}\}\le7.$$
In conclusion,
\begin{align*}
\sharp(C(\bar K)\cap\Gamma)\le&(16g^2+32g+124)(20g)^\rho+14\cdot13^{\rho+1}\\
\le&(16g^2+32g+184)(20g)^\rho.
\end{align*}
This proves the first assertion. The second one is the special case that $\Gamma=J(K)$.
\end{proof}

\bibliographystyle{alphaurl}

\bibliography{bib}

\end{document}